\documentclass[a4paper,12pt]{article}
\usepackage{pstricks}
\usepackage{pstricks-add}
\usepackage{amsmath, amsfonts, amssymb, amstext, amsthm} 
\usepackage{graphicx, subfigure} 
\usepackage[affil-it]{authblk} 
\usepackage{tikz}

\def\R{\mathbb{R}}

\def\e{\varepsilon}

\def\la{\lambda}

\def\nub{\overline{\nu}}
\def\mub{\overline{\mu}}

\def\fp{f'(0)}

\def\vp{\varphi}
\def\Gt2{\widetilde{\Gamma_2}}

\def\lp{\left(}
\def\rp{\right)}

\def\lV{\left\Vert}
\def\rV{\right\Vert}
\def\lc{\left[}
\def\rc{\right]}

\newtheorem{theorem}{Theorem}[section]
\newtheorem{prop}{Proposition}[section]

\theoremstyle{definition}
\newtheorem{definition}{Definition}[section]

  \author{Antoine Pauthier%
  \thanks{e-mail: \texttt{antoine.pauthier@math.univ-toulouse.fr}}}
\affil{Institut de Math\'ematiques de Toulouse ; UMR5219 \\ Universit\'e de Toulouse ; CNRS \\ UPS IMT, F-31062 Toulouse Cedex 9, France}

\title{Road-field reaction-diffusion system: a new threshold for long range exchanges}

\begin{document}
\maketitle

\begin{abstract}
We consider reaction-diffusion equations of KPP type in a presence of a line of fast diffusion with non-local exchange terms
between the line and the framework. Our study deals with the infimum of the spreading speed depending on the exchange functions.
We exhibit a new threshold in the limit of long range exchange terms for the line to influence the propagation.
\end{abstract}

\section{Introduction}
The purpose of this note is to study some properties concerning the spreading speed of the 
following reaction-diffusion system, introduced in \cite{Pauthier}:
\begin{equation}
 \label{RPeq}
\begin{cases}
 \partial_t u-D \partial_{xx} u = -\mub u+\int \nu(y)v(t,x,y)dy & x \in \R,\ t>0 \\
 \partial_t v-d\Delta v = f(v) +\mu(y)u(t,x)-\nu(y)v(t,x,y) & (x,y)\in \R^2,\ t>0.
\end{cases}
\end{equation}
The initial road-field system was introduced in \cite{BRR1}. It was generalised to nonlocal exchange terms in \cite{Pauthier} and \cite{Pauthier2}.
We refer to these papers for more informations. We use the notation
$\mub=\int\mu,$ $\nub=\int\nu.$ Thus, it is easy to check that without reaction, the above system is mass-conservative.
Our assumptions are the following.
\begin{itemize}
 \item The reaction term $f$ is of KPP type, i.e. strictly concave with $f(0)=f(1)=0,$ and quadratic outside $[0,1].$
 \item The two exchange functions $\mu$ and $\nu$ are continuous, nonnegative, even. For the sake of simplicity, we will consider compactly supporded functions,
 but our results can easily be extended to a mere general class of functions. See \cite{Pauthier} for the optimal (to our knowldge) hypothesis.
\end{itemize}
The purpose of the model (\ref{RPeq}) is to study a propagation driven by the line.
This is the main motivation of the following Theorem, which also gives a definition of the \textbf{spreading speed} 
for this kind of model. It was proved in \cite{Pauthier}.

 \begin{theorem}\label{spreadingthmgeneral}
Let $(u,v)$ be a solution of (\ref{RPeq})
with a nonnegative, compactly supported initial datum $(u_0,v_0)$.
Then, there exists an asymptotic speed of spreading $c^*$ 
and a unique positive bounded stationary solution of (\ref{RPeq})
$(U,V)$ such that, 
pointwise in $y$, we have: 
 \begin{itemize}
  \item for all $c>c^*$, $\lim_{t\to\infty}\sup_{|x|\geq ct}(u(t,x),v(t,x,y)) = (0,0)$ ;
  \item for all $c<c^*$, $\lim_{t\to\infty}\inf_{|x|\leq ct}(u(t,x),v(t,x,y)) = (U,V)$.
 \end{itemize}
\end{theorem}

\section{Infimum for the spreading speed}
 For fixed parameters $d,D,\fp,\mub,\nub$ we consider the set of admissible exchanges
$$
\Lambda_{\mub} = \{\mu\in C_0(\R),\mu\geq 0, \int\mu=\mub,\mu \textrm{ even} \}.
$$
We define $\Lambda_{\nub}$ in a similar fashion.
For $\mu\in\Lambda_{\mub}$ and $\nu\in\Lambda_{\nub},$ there exists a spreading speed $c^*(\mu,\nu).$ 
A natural question is to wonder about the existence of maximal or minimal spreading speed for $\mu,\nu$ admissible exchanges. 
This note is devoted to the existence of an infimum for the spreading speed. Thus, we also prove that there is no minimal spreading speed.
The main result relies on the following theorem.

\begin{theorem}\label{ralentissement}
Let us consider the nonlocal system (\ref{RPeq}) with fixed exchange masses $\mub$ and $\nub.$
Let $c^*$ be the spreading speed given by Theorem \ref{spreadingthmgeneral}, depending on the repartition of $\mu$ or $\nu.$
	\begin{enumerate}
	\item If $\displaystyle D\in \lc 2d,d\lp 2+\frac{\mub}{\fp}\rp\rc,$ $\displaystyle \inf c^*=2\sqrt{d\fp}.$
	\item Fix $\displaystyle D>d\lp 2+\frac{\mub}{\fp}\rp,$ then $\displaystyle \inf c^*>2\sqrt{d\fp}.$
	\end{enumerate}
Moreover, in both cases, minimizing sequences can be given by long range exchange terms of the form 
$\mu_R(y)=\frac{1}{R}\mu\lp\frac{y}{R}\rp$ or $\nu_R(y)=\frac{1}{R}\nu\lp\frac{y}{R}\rp$ with $R\to \infty.$
\end{theorem}

Let us recall that in \cite{BRR1} and \cite{Pauthier} was exhibited the threshold $D=2d$ for the spreading in
the direction of the road, whatever be the considered road-field system:
\begin{itemize}
\item if $D\leq 2d,$  $c^*=c_K:=2\sqrt{d\fp}$;
\item if $D>2d,$ $c^*>c_K.$
\end{itemize}
We show that the threshold $\displaystyle D=2d+d\frac{\mub}{\fp}$ has an important effect.
\begin{itemize}
 \item If $\displaystyle D>2d+d\frac{\mub}{\fp},$ the speed $c^*$ is strictly greater than $c_K$ in the $x$-direction,
 with a bound independent of the exchange functions.
 \item If $\displaystyle D<2d+d\frac{\mub}{\fp},$ $c^*$ tends to $c_K$ as the exchange functions vanish.
\end{itemize}

\section{Background on the computation of the spreading speed}
The importance of the linearised system for KPP-type model motivates the following definition for traveling waves.
\begin{definition}\label{deftravelling}
We call a \textbf{linear traveling wave} a 3-tuple
 $(c,\la,\phi)$ with $c>0,$ $\la>0,$ and $\phi\in H^1(\R)$ a positive function such that
 $$
 \begin{pmatrix}
  u \\ v
 \end{pmatrix}
\mapsto e^{-\la(x-ct)} \begin{pmatrix}
                        1 \\ \phi(y)
                       \end{pmatrix}
$$
 be a solution of the corresponding linearised system in 0. $c$ is the speed of the exponential traveling waves.
\end{definition}
The previous definition for traveling waves provides us a helpful characterisation for spreading speed.
\begin{prop}\label{defspreadingspeed}
The \textbf{spreading speed $c^*$} given by Theorem \ref{spreadingthmgeneral}
can be defined as follows:
$$
c^*=\inf\{c>0 | \text{ linear traveling waves with speed $c$ exists}\}.
$$
\end{prop}

Inserting definition supplied by Proposition \ref{defspreadingspeed} into (\ref{RPeq}) yields the following system 
in $(c,\la,\phi):$
\begin{equation}
\label{systemlambdaphi}
\begin{cases}
 -D\la^2+\la c+\mub  = \int \nu(y)\phi(y)dy \\
-d\phi''(y)+(\la c-d\la^2-f'(0)+\nu(y))\phi(y)  =  \mu(y).
\end{cases}
\end{equation}
These equations and integrals have to be understood in a distribution sense if needed.
As explained in \cite{Pauthier}, the first equation of (\ref{systemlambdaphi}) 
gives the graph of a function $\la\mapsto\Psi_1(\la,c):= -D\la^2+\la c+\mub$, which
means to be equal to $\int \nu(y)\phi(y)dy$, provided $(c,\la,\phi)$ defines an exponential traveling waves. 

The second equation of (\ref{systemlambdaphi})
gives, under some assumptions on $\la$, a unique nonnegative solution $\phi=\phi(y;\la,c)$ in $H^1(\R)$. To this unique solution
we associate the function $\Psi_2(\la,c):=\int \nu(y)\phi(y)dy$.
Let us denote $\Gamma_1$ the graph of $\Psi_1$ in the $(\la, \Psi_1(\la))$ plane, and $\Gamma_2$ the graph
of $\Psi_2$. So, (\ref{systemlambdaphi}) amounts to the investigation of $\la,c>0$ such that $\Gamma_1$ and $\Gamma_2$ intersect.

\subsection{Resolution of the $(c,\la,\phi)$-system: general remarks}
Thereafter we recall some facts on the two functions $\Psi_{1,2}.$
For more details and proofs, we refer to \cite{Pauthier}. 
\paragraph{Behaviour of $\Psi_1$}
Let us recall that in a $(\la,\Psi_1(\la))$ plane, $\Psi_1$ defines parabola, nonnegative
for $\displaystyle \la\in[\la_1^-,\la_1^+]$ with $\displaystyle \la_1^\pm = \frac{c\pm\sqrt{c^2+4D\mub}}{2D},$ and that 
$$\Psi_1(0)=\Psi_1(\frac{c}{D})=\mub.$$
Notice that $\Psi_1$ depends on $c,D,$ and $\mub,$ but does not depend on $\nu,d,$ neither on the repartition of $\mu.$
Only the mass matters. This will be useful for the sequel.


\paragraph{Behaviour of $\Psi_2$}
The function $\Psi_2$ is defined implicitly by the solution of the following system (\ref{eqgeneralesurphi})
\begin{equation}
\label{eqgeneralesurphi}
\begin{cases}
 -d\phi''(y)+(\la c-d\la^2-f'(0)+\nu(y))\phi(y)  =  \mu(y) \\
 \phi \in H^1(\R),\ \phi\geq0. 
\end{cases}
\end{equation}
If it exists, to the solution of (\ref{eqgeneralesurphi}) we associate $\displaystyle\Psi_2(\la):=\int_\R \nu(y)\phi(y)dy.$
It has been shown that $\Psi_2$ is defined for $c>c_K$ and $\displaystyle\la\in]\la_2^-,\la_2^+[,$ with
$$
\la_2^\pm=\frac{c\pm\sqrt{c^2-c^2_K}}{2d}.
$$
Recall that the classical Fisher-KPP speed is given by $c_K=2\sqrt{d\fp}.$
$\Psi_2$ is a smooth convex function, symmetric with respect to the line $\{\la=\frac{c}{2d}\},$ and can be continuously extended
to $\la_2^\pm$ by $\Psi_2(\la_2^\pm)=\mub,$ with vertical tangents at these points.

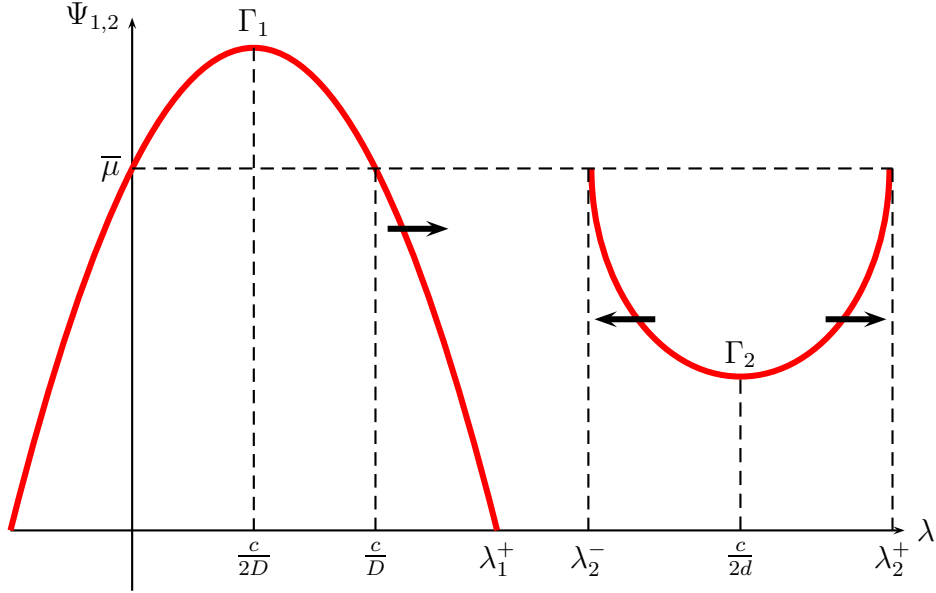
\begin{figure}[!ht]
   \centering
   \psset{unit=0.8}
\def\mub{\overline{\mu}}

\begin{pspicture}(-3,-1)(13,9)
\psline{->}(-2,0)(12.7,0)
\psline{->}(0,-1)(0,8.5)
\psparabola[linewidth=0.1,linecolor=red](-2,0)(2,8)
\psellipticarc[linewidth=0.1,linecolor=red](10,6)(2.5,3.5){180}{360}

\psline[linestyle=dashed](0,6)(12.5,6)
\uput{0.2}[180](0,6){$\mub$}
\uput{0.2}[0](12.7,0){$\lambda$}
\uput{0.2}[-90](6,0){$\lambda_1^+$}
\uput{0.2}[-90](7.5,0){$\lambda_2^-$}
\uput{0.2}[-90](12.5,0){$\lambda_2^+$}
\psline[linestyle=dashed](2,8)(2,0)
\psline[linestyle=dashed](4,6)(4,0)
\psline[linestyle=dashed](7.5,6)(7.5,0)
\psline[linestyle=dashed](10,2.5)(10,0)
\psline[linestyle=dashed](12.5,6)(12.5,0)
\uput{0.2}[-90](10,0){$\frac{c}{2d}$}
\uput{0.2}[-90](2,0){$\frac{c}{2D}$}
\uput{0.2}[-90](4,0){$\frac{c}{D}$}
\uput{0.2}[90](10,2.5){$\Gamma_2$}
\uput{0.2}[90](2,8){$\Gamma_1$}
\uput{0.2}[180](0,8.5){$\Psi_{1,2}$}
\psline[linewidth=0.1]{->}(4.2,5)(5.2,5)
\psline[linewidth=0.1]{->}(8.6,3.5)(7.6,3.5)
\psline[linewidth=0.1]{->}(11.4,3.5)(12.4,3.5)
\end{pspicture}
   \caption{\label{2graphs}representation of $\Gamma_1$ and $\Gamma_2,$ behaviours as $c$ increases}
\end{figure}

Notice that, contrary to $\Psi_1,$ $\Psi_2$ is highly dependent on the two exchange functions $\nu$ and $\mu.$
Thus, this paper will mainly focus on this function and how its variations depend on these exchange functions.
However, the extreme points $(\la_2^\pm,\Psi_2(\la_2^\pm))$ do not depend on these functions, but only on $c,d,\fp,$ and $\mub.$
Global behaviours of $\Psi_1$ and $\Psi_2$ are summarised in Figure \ref{2graphs}.

\paragraph{Intersection of $\Gamma_1$and $\Gamma_2$}  We focus on the case $D>2d,$ the other one leading to dynamics with no influence of the road.
The functions of the $c$-variable $c\mapsto\Psi_1$ and 
$c\mapsto\la_2^-$ are respectively increasing and decreasing. Hence, there exists $c^*$
such that $∀c > c^*,$ $\Gamma_1$ and $\Gamma_2$ intersect, and $∀c < c^*,$ $\Gamma_2$ does not intersect the closed
convex hull of $\Gamma_1.$ Moreover, the strict concavity of $\Gamma_1$ and the strict convexity of $\Gamma_2$ 
give that for $c = c^*,$ $\Gamma_1$ and $\Gamma_2$ are tangent on $\la^*$ and for $c > c^*,$ $c$ close to
$c^*,$ $\Gamma_1$ and $\Gamma_2$ intersect twice.

As the extreme points $(\la_2^\pm,\Psi_2(\la_2^\pm))$ are independent of $\nu$ and the repartition of $\mu$, the study of the spreading speed associated to two
exchange functions amounts to analyse the relative position of the corresponding $\Psi_2$ functions. Given the monotonicity in $c,$ 
a $\Psi_2$ function under another leads to a slower spreading speed, and vice versa. The convexity of $\Psi_2$ even allows us to
study the variations in a neighbourhood of $(\la,c^*)$ to get a local result.

\section{A new threshold}
We show how long range exchange terms tend to slow down the dynamics.
More precisely, for any given functions $\mu,\nu,$ we set
\begin{equation}\label{longrangeExc}
\mu_R(y)=\frac{1}{R}\mu\lp\frac{y}{R}\rp,\ \nu_R(y)=\frac{1}{R}\nu\lp\frac{y}{R}\rp.	
\end{equation}
This asymptotics yields to a new threshold in order to get or be greater than the KPP spreading speed.
Moreover, this provides minimizing sequences for the spreading speed, as asserted in Theorem \ref{ralentissement}.

In the system (\ref{RPeq}), replace at least one exchange function by a long range exchange function given by (\ref{longrangeExc}),
and let us denote $c^*(R)$ the corresponding spreading speed in the sense of Proposition \ref{defspreadingspeed}.


The proof relies on a very simple remark: for $\Gamma_1$ and $\Gamma_2$ to intersect, it is necessary to have $\la_1^+>\la_2^-$ 
(see Figure \ref{2graphs}).
Let us see $\la_1^+$ and $\la_2^-$ as functions of the speed $c$ given by
\begin{align}
	\la_1^+ : c \longmapsto & \frac{c+\sqrt{c^2+4D\mub}}{2D}  \label{lambda1} \\
	\la_2^- : c  \longmapsto & \frac{c-\sqrt{c^2-c_K^2}}{2d}.  \label{lambda2}
\end{align}
They are both continuous. $\la_1^+$ is increasing on $[0,+\infty[,$ $\la_2^-$ is decreasing on $[c_K,+\infty[.$ We may also notice that 
$\la_1^+$ is decreasing with respect to $D.$ An
explicit computation gives
\begin{equation}\label{threshold}
	d\lp 2+\frac{\mub}{\fp}\rp=\inf \{D>0,\ \la_1^+(c_K)<\la_2^-(c_K)\}.
\end{equation}

\begin{figure}[!ht]
   \centering
   \psset{unit=0.9}
\def\mub{\overline{\mu}}

\begin{pspicture}(-3,-1)(13,9)
\psline{->}(-2,0)(11.3,0)
\psline{->}(0,-1)(0,8.5)
\psparabola[linewidth=0.1,linecolor=blue](-2,0)(2,8)
\psellipticarc[linewidth=0.1,linecolor=red](8.5,6)(2.5,3.5){180}{360}

\psline[linestyle=dashed](0,6)(11,6)
\uput{0.2}[180](0,6){$\mub$}
\uput{0.2}[0](11.1,0){$\lambda$}
\uput{0.2}[-90](6,0){$\lambda_1^+=\la_2^-$}
\psline[linestyle=dashed](2,8)(2,0)
\uput{0.2}[-90](2,0){$\frac{c}{2D}$}
\uput{0.2}[90](8.5,2.5){$\Gamma_2$}
\uput{0.2}[90](2,8){$\Gamma_1$}
\psline[linewidth=0.1]{->}(8.5,2.5)(8.5,1.5)

\psline[linewidth=0.1,linecolor=red](6,6)(6,0)
\psline[linewidth=0.1,linecolor=red](6,0)(11,0)
\psline[linewidth=0.1,linecolor=red](11,0)(11,6)
\end{pspicture}
   \caption{\label{Rinfini}Behaviour of $\Gamma_2$ as $R\to+\infty,$ critical case $\la_1^+=\la_2^-$}
\end{figure}
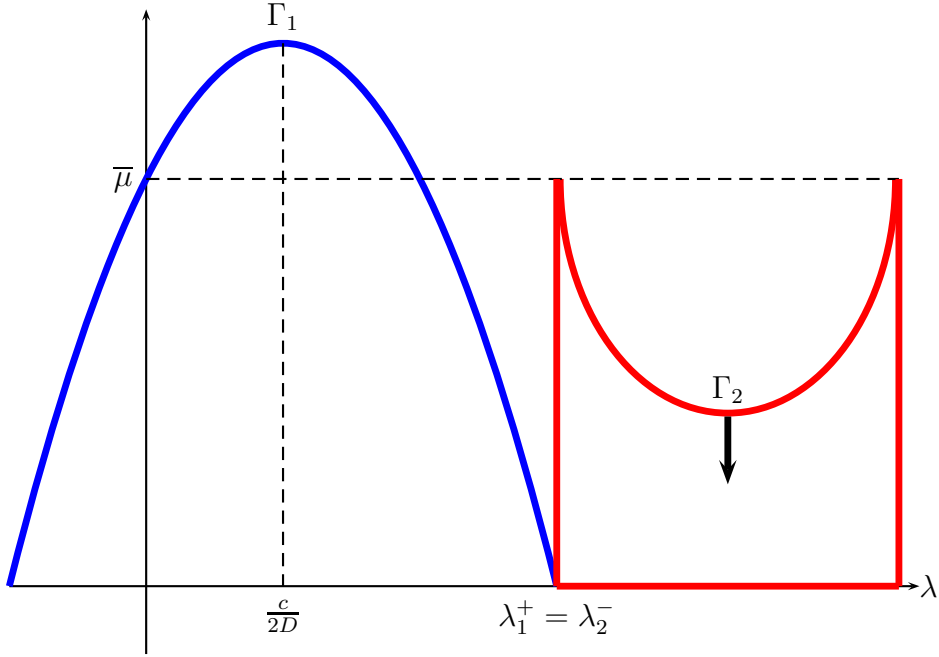

\subparagraph{Proof of the first part of Theorem \ref{ralentissement}} For the sake of simplicity we will focus on the
general model (\ref{RPeq}), the other being similar and even easier - see \cite{Pauthier}. 
Let  $D$ be less than $\displaystyle d\lp 2+\frac{\mub}{\fp}\rp,$ $\e>0.$ Let $c\in]c_K,c_K+\e[.$ Then
$\la_1^+(c)>\la_2^-(c).$ Choose any $\la_0\in]\la_2^-,\la_1^+[.$

\textit{First case: long range for $\mu.$} Let us replace $\mu$ by $\displaystyle\mu_R(y)=\frac{1}{R}\mu\lp\frac{y}{R}\rp$ in (\ref{RPeq}).
The $(c,\la,\phi)$  associated equation (\ref{eqgeneralesurphi}) is
\begin{equation}
\label{eqphimuR}
\begin{cases}
-d\phi_R''(y)+(P(\la)+\nu(y))\phi_R(y)  =  \frac{1}{R}\mu\lp\frac{y}{R}\rp\\
\phi \in H^1(\R),\ \phi\geq0
\end{cases}
\end{equation}
where as usual $P(\la)=\la c-d\la^2-f'(0).$ The curve $\Gamma_2$ is defined as the graph of 
$$\Psi_2:\la \longmapsto \int_\R\nu\phi_R$$
where $\phi_R$ is the unique solution of (\ref{eqphimuR}). From the choice of $\la_0,$ $P(\la_0)>0.$ The maximum principle yields
for $\la=\la_0$
$$
\lV\phi_R\rV_\infty \leq \frac{1}{RdP(\la_0)}\lV\mu\rV_\infty 
$$
which gives
\begin{equation*}
\Psi_2(\la_0)\leq \frac{\nub}{RdP(\la_0)}\lV\mu\rV_\infty \underset{R\to+\infty}{\longrightarrow}0.
\end{equation*}
Thus, there exists $R_0,$ $\forall R>R_0,$ $\Psi_1(\la_0)>\Psi_2(\la_0),$ so $c^*(R)<c.$

\textit{Second case: long range for $\nu.$} The study is quite similar. The $(c,\la,\phi)$  associated equation (\ref{eqgeneralesurphi}) is
\begin{equation}
	\label{eqphinuR}
	\begin{cases}
		-d\phi_R''(y)+(P(\la)+\frac{1}{R}\nu\lp\frac{y}{R}\rp)\phi_R(y)  =  \mu(y) \\
		\phi \in H^1(\R),\ \phi\geq0
	\end{cases}
\end{equation}
Let $\vp$ be the only $H^1$ solution of 
$$
-d\vp''(y)+P(\la_0)\vp(y)  =  \mu(y).
$$
This provides a supersolution for (\ref{eqphinuR}) with $\la=\la_0.$ As $\mu$ is compactly supported, $\vp$ belongs to $L^1(\R).$ Hence
\begin{equation*}
\Psi_2(\la_0)=\frac{1}{R}\int_{\R}\nu(\frac{y}{R})\phi_R(y)dy\leq\frac{\lV\nu\rV_\infty}{R}\int_{\R}\vp(y)dy\underset{R\to+\infty}{\longrightarrow}0
\end{equation*}
and we conclude as in the previous case.

\subparagraph{Proof of the second part of Theorem \ref{ralentissement}}
Let  $D$ be greater than $\displaystyle d\lp 2+\frac{\mub}{\fp}\rp.$ 
Thus, $\la_1^+(c_K)<\la_2^-(c_K).$ 
From the above reasoning, the minimal speed $c_{min}$ is given by
$$
c_{min} =\inf\{c,\ \la_1^+(c)>\la_2^-(c)\}.
$$
Continuity and monotonicity of $\la_1^+,\la_2^-$ given by (\ref{lambda1})-(\ref{lambda2}) 
warrant the existence of $c_{min}$ and the inequality $c_{min}>c_{K}.$ 
Moreover, the above study ensures us that it is optimal. See Figure \ref{Rinfini}.
\qed

\section{Remarks and opened questions}
The above result can easily be extended to the semi-limit models presented in \cite{Pauthier}, with one nonlocal exchange and the other exchange 
by boundary condition.

Using the same kind of geometric considerations, it is easy to give the following upper bound for the spreading speed.
\begin{prop}
 For fixed parameters $d,D,\mub,\nub,\fp,$ then for all admissible exchanges $\mu\in\Lambda_{\mub}$ and $\nu\in\Lambda_{\nub}$ we 
 have
 $$
 c^*(\mu,\nu)\leq D\sqrt{\frac{\fp}{D-d}}.
 $$
\end{prop}
An opened question is to know if this bound is reached for some exchanges. Moreover, 
this bound ensures us the existence of minimizing sequences for the spreading speed. 
Hence, it is questionable whether these sequences converge, in which sense, etc.

\paragraph{Acknowledgements} The research leading to these results has received 
funding from the European Research Council under the European Union’s 
Seventh Framework Programme (FP/2007-2013) / ERC Grant Agreement n.321186 - ReaDi -Reaction-Diffusion Equations, Propagation and Modelling.
I am also grateful to Henri Berestycki and Jean-Michel Roquejoffre for suggesting me the
models and many fruitful conversations.

\bibliographystyle{plain}
\footnotesize
\bibliography{biblio} 
 
\end{document}